\documentclass[11pt]{amsart}
   \usepackage{amsmath,amssymb}
   
   \renewcommand{\bf}{\bfseries}

   \newcommand{\E}{\mathbb{E}}

   \newtheorem{theorem}{Theorem}
   \newtheorem{lemma}{Lemma}

   \renewcommand{\epsilon}{\varepsilon}

   \newcommand{\dis}{\displaystyle}
  
\begin{document}
  \title{Pointwise convergence of averages along cubes II}
   \author{I. Assani}
 \begin{abstract}
   Let $(X,\mathcal{B},\mu, T)$ be a measure preserving system.
   We prove the pointwise convergence of averages along cubes of
   $2^{k}-1$ bounded and measurable functions for all $k$.

  \end{abstract}
  \maketitle
  \section{Introduction}

 Let $(X, \mathcal{B}, \mu, T)$ be a dynamical system where T is a measure preserving transformation on the measure space
 $(X, \mathcal{B}, \mu, T)$. In \cite{assani} we proved the pointwise convergence of the averages
  $$\frac{1}{N^2}\sum_{n,m= 0}^{N-1}
    f_1(T^nx)f_2(T^mx)f_3(T^{n+m}x)$$ and of similar averages
    with seven bounded functions $f_i$. We also showed that if $T$ is weakly mixing then similar averages for $2^{k}-1$ bounded functions
converge a.e to the product of the integrals of the functions
$f_i$.
 The averages of three functions were used  in
\cite{Bergelson} to generalize Khintchine recurrence result
\cite{Khintchine}. In \cite{Host-Kra} B.Host and B.Kra proved that
the averages of $2^k -1$ bounded functions converge in $L^2$ norm.
 To achieve this result they identified increasing factors $Z_k$, $k=0,
1,2,...$ of ergodic dynamical systems and showed the following
\begin{itemize}
\item The averages of $2^k -1$ bounded functions converge a.e if
each function belongs to the factor $Z_{k-1}$. They used for that
a result of A. Leibman \cite{Leibman}.
 \item The averages of $2^k -1$ functions converge in $L^2$ norm
\end{itemize}
One consequence of their method is that for each $k$ the factor
$Z_{k-1}$ is characteristic for the $L^2$ norm of the averages of
$2^{k-1}$ functions. Let us note that $Z_1$ is the Kronecker
factor and $Z_2$ the CL factor. The notion of characteristic
factor is due to H. Furstenberg and can be found explicitly stated
in \cite{Furstenberg-Weiss}. Our main results are the following

\begin{theorem}
 Let $(X, \mathcal{B}, \mu, T)$ be a measure preserving system .
 Then the averages of the cubes of $2^k -1$ functions converge
 a.e.
 \end{theorem}
  One consequence of the path we use is the following

  \begin{theorem}
  Let $(X, \mathcal{B}, \mu, T)$ be an ergodic dynamical system.
  For each $k\geq 1$ the factors $Z_{k-1}$ is characteristic for the
  pointwise convergence of the averages along the cubes of $2^k-1$ bounded
  and measurable functions.
  \end{theorem}
  \section{Inequalities from the averages of three or seven functions}

   As shown in \cite{Host-Kra} the factors $Z_k$ can be defined inductively by using the
   seminorms $|||.|||_k$ where
   \begin{enumerate}
   \item $|||f|||_1 = \big|\int fd\mu\big|$
   \item For every $k\geq 1$
   $$|||f|||_{k+1}^{2^{k+1}} = \lim_H
   \frac{1}{H}\sum_{h=1}^{H}|||f.f\circ T^h|||_{k}^{2^{k}}$$
   \end{enumerate}

   With the help of these semi norms factors are built with the property that
   for all $f\in L^{\infty}$ we have $\E(f|Z_k) = 0$ if and only
   if $|||f|||_{k+1}= 0$.
   \smallskip

    We mention a few inequalities that were used in \cite{assani} in the proof of
    the pointwise convergence of the averages of three and seven
    bounded functions. The constant $C$  may change
 from one line to the other. But it will depend only at time on the $L^{\infty}$ norm of the functions $f_j$.
 For all bounded functions $f_i$, $1\leq i\leq 7$ we have

\begin{equation}
  \bigg(\frac{1}{N}\sum_{n=0}^{N-1}\big|\frac{1}{N}\sum_{m=0}^{N-1}f_2(T^mx)f_3(T^{n+m}x)\big|^2\bigg)\leq
  C\sup_t\bigg|\frac{1}{N}\sum_{m'=0}^{N-1}f_3(T^{m'}x)e^{2\pi
im't}\bigg|^2\|f_2\|_{\infty}^2.
  \end{equation}
\begin{equation}
  \limsup_N \sup_t \big|\frac{1}{N}\sum_{n=0}^{N-1}f(T^nx)e^{2\pi
  int}\big|^2 \leq C\||f|\|_2^2.
  \end{equation}
\begin{equation}
\begin{aligned}
&\frac{1}{N^2}\sum_{p=0}^{N-1}\sum_{n=0}^{N-1}\|f_1\|_{\infty}^2\|f_2\|_{\infty}^2\|f_3\|_{\infty}^2\bigg|\frac{1}{N}
\sum_{m=0}^{N-1}f_4(T^mx)f_5(T^{n+m}x)f_6(T^{p+m}x)f_7(T^{p+n+m}x)\bigg|^2\\
&\leq C\prod_{i=1}^{5}\|f_i\|_{\infty}^2
\frac{1}{N}\sum_{n=0}^{N-1}\sup_t
\bigg|\frac{1}{N}\sum_{m'=0}^{N-1}f_6(T^{m'}x)f_7(T^{n+m'}x)e^{2\pi
im't}\bigg|^2
\end{aligned}
\end{equation}
\begin{equation}
\limsup_N\frac{1}{N}\sum_{n=0}^{N-1}\sup_t\bigg|\frac{1}{N}\sum_{m=0}^{N-1}f_1(T^mx)f_2(T^{n+m}x)e^{2\pi
imt}\bigg|^2\leq C Min[\||f_1|\|_3^2, \||f_2|\|_3^2].
\end{equation}
 With these inequalities we can explain the first induction step
 allowing to get the convergence of the averages of seven
 functions from the inequalities obtained for the averages of
 three functions.  We hope that these explanations will make the
 proof of theorem 1 more transparent.
 The square of the averages of three functions $$M_N(f_1, f_2, f_3)(x)= \frac{1}{N^2}\sum_{n,m= 0}^{N-1}
    f_1(T^nx)f_2(T^mx)f_3(T^{n+m}x)$$ is bounded
by
$$C\bigg(\frac{1}{N}\sum_{n=0}^{N-1}\big|\frac{1}{N}\sum_{m=0}^{N-1}f_2(T^mx)f_3(T^{n+m}x)\big|^2\bigg).$$
which by the equation (1) is bounded by
$$C\sup_t\bigg|\frac{1}{N}\sum_{m'=0}^{N-1}f_3(T^{m'}x)e^{2\pi
im't}\bigg|^2\|f_2\|_{\infty}^2.$$ The equation (2) guarantees
that the $\limsup $ of this last quantity is equal to zero if
$$\||f_3|\|_2 = 0$$. This is the same of saying that the function
$f_3\in \mathcal{K}^{\perp}.$ By using a similar path for the
functions $f_1$ and $f_2$ one can see that the Kronecker factor is
pointwise characteristic for the averages of three bounded
functions.
\smallskip

The square of the averages of seven bounded functions $f_i$,
$1\leq i\leq 7$, is bounded by
$$
\frac{1}{N^2}\sum_{p=0}^{N-1}\sum_{n=0}^{N-1}\|f_1\|_{\infty}^2\|f_2\|_{\infty}^2\|f_3\|_{\infty}^2\bigg|\frac{1}{N}
\sum_{m=0}^{N-1}f_4(T^mx)f_5(T^{n+m}x)f_6(T^{p+m}x)f_7(T^{p+n+m}x)\bigg|^2.
$$
By using the equation (3) this term is bounded by
$$\leq C\prod_{i=1}^{5}\|f_i\|_{\infty}^2
\frac{1}{N}\sum_{n=0}^{N-1}\sup_t
\bigg|\frac{1}{N}\sum_{m'=0}^{N-1}f_6(T^{m'}x)f_7(T^{n+m'}x)e^{2\pi
im't}\bigg|^2
$$
Then the equation (4) applied to $f_6$ and $f_7$ shows that
$$\limsup_N \frac{1}{N}\sum_{n=0}^{N-1}\sup_t
\bigg|\frac{1}{N}\sum_{m'=0}^{N-1}f_6(T^{m'}x)f_7(T^{n+m'}x)e^{2\pi
im't}\bigg|^2= 0,
$$
if one of the functions $f_6$ or $f_7$ belongs to $CL^{\perp}$.
The equation (4) is a consequence of the equation (1) and of the
van der Corput's inequality. In fact it is by averaging with h
means of functions of the form $f.f\circ T^h$ that the semi norms
$\||f|\|_3$ appear.
\smallskip
We are going to follow a similar path to prove theorem 1.
\section{Proof of theorem 1}
  We will prove theorem 1 by induction on k. In \cite{assani} we
  proved that the averages of seven functions converge a.e. We
  showed that the $Z_2 = CL$ factor was characteristic for the
  pointwise convergence of seven functions. This establishes the first step of the induction process.
  We will use the same
  notation and some of the remarks made in \cite{assani}.
  \begin{itemize}
  \item For each $k\geq 4$ we denote by
  $M_N(f_1, f_2, ...,f_{2^k-1})$ the averages of $2^k-1$ bounded
  functions. Without loss of generality we assume that the
  functions are bounded by 1 in absolute value.
 \item The functions $f_j$ are listed in such a way that those
 depending on the index $i_k$ are indexed by those $j$ ,
 $2^{k-1}\leq j \leq 2^k - 1$. The product of these terms
 depending on $i_k$ is denoted by $S_{N, (i_1, i_2,...,
 i_k)}(f_{2^{k-1}}, ..., f_{2^k -1})(x)$. Each term
$S_{N, (i_1, i_2,...,
 i_k)}(f_{2^{k-1}}, ..., f_{2^k -1})(x)$ is the product of two
 groups of $2^{k-2}$ functions denoted by
 $$A_{N, (i_1, i_2,..., i_k)}(f_{2^{k-1}}, f_{2^{k-1}
 +1},...,f_{3.2^{k-2}}(x)$$ and $$B_{N, (i_1, i_2,...,
 i_k)}(f_{3.2^{k-2} +1},..., f_{2^k -1})(x)$$ where the powers of $T$
 associated with each function in the second group are those
 appearing in the first group shifted by the index $i_1$.
 We have
 $$B_{N, (i_1, i_2,...,
 i_k)}(f_{3.2^{k-2} +1},..., f_{2^k -1})(x)= A_{N, (i_1, i_2,...,
 i_k)}(f_{3.2^{k-2}+1},...,f_{2^k-1})(T^{i_1}x)$$
 \item We have also the inequality
 \begin{equation}
\begin{aligned}
&|M_N(f_1, f_2,..., f_{2^k -1})(x)|^2 \\
&\leq \prod_{j=1}^{2^{k-1}-1}\|f_j\|_{\infty}^2
\frac{1}{N^{k-1}}\sum_{i_1,...,
 i_{k-1}=0}^{N-1}\big|\frac{1}{N}\sum_{i_k=0}^{N-1}S_{N,(i_1,i_2,...,i_k)}(f_{2^{k-1}},..., f_{2^k
-1})(x)\big|^2.
\end{aligned}
\end{equation}
 \end{itemize}

\noindent{\bf Induction Assumption} \vskip1ex We make the
following assumption
 \vskip1ex For all bounded functions $g_j$,
$3.2^{k-2}+1 \leq j \leq 2^{k}-1$ we have
\begin{equation}
\begin{aligned}
&\limsup_N\frac{1}{N^{k-2}}\sum_{i_1,..., i_{k-2}= 0}^{N-1}
\bigg|\frac{1}{N}\sum_{i_k=0}^{N-1}A_{N,(i_1,i_2,...,i_{k-2},i_k)}(g_{3.2^{k-2}+1}.,...,g_{2^k
-1})(x)\bigg|^2 \\
&\leq C.Min_{\{3.2^{k-2}+1 \leq j \leq
2^{k}-1\}}\||g_j|\|_{k-2}^2.
\end{aligned}
\end{equation}

 As indicated above this assumption is shown to be true for $k = 3, 4$ in
 \cite{assani}. We want to show that it also holds for k. To this
 end we have the following extension of lemma 4 in \cite{assani}.
 \begin{lemma}
 If one of the $2^{k-2}$ functions $f_j$, $3.2^{k-2}+1\leq j\leq
2^{k}-1$ is in $Z_{k-1}^{\perp}$ then
\begin{equation}
\lim_N\frac{1}{N^{k-2}}\sum_{i_1,..., i_{k-2}=
0}^{N-1}\sup_{t}\bigg|\frac{1}{N}\sum_{i_k=0}^{N-1}A_{N,
(i_1,i_2,...,i_{k-2},i_k)}(f_{3.2^{k-2}+1},...,f_{2^k
-1})(x)e^{2\pi ii_kt}\bigg|^2 =0
\end{equation}
\end{lemma}
\begin{proof}
 We use now the same path as in \cite{assani}. With Van der Corput lemma applied to
 each term
 \[
 \sup_{t}\bigg|\frac{1}{N}\sum_{i_k=0}^{N-1}A_{N, (i_1,i_2,...,i_{k-2},i_k)}(f_{3.2^{k-2}+1},...,f_{2^k
-1})(x)e^{2\pi ii_kt}\bigg|^2,
\]
we have then for each $(H + 1)<< N$
\[
\begin{aligned}
&\frac{1}{N^{k-2}}\sum_{i_1,..., i_{k-2}=
0}^{N-1}\sup_{t}\bigg|\frac{1}{N}\sum_{i_k=0}^{N-1}A_{N,(i_1,i_2,...,i_{k-2},i_k)}(f_{3.2^{k-2}+1},...,f_{2^k
-1})(x)e^{2\pi ii_kt}\bigg|^2 \\
&\leq C.\bigg(\frac{1}{H} + \frac{1}{H}\sum_{h=1}^H
\frac{1}{N^{k-2}}\sum_{i_1,..., i_{k-2}=
0}^{N-1}\\
&\bigg|\frac{1}{N}\sum_{i_k=0}^{N-h-1}A_{N,(i_1,i_2,...,i_{k-2},i_k)}(f_{3.2^{k-2}+1}.f_{3.2^{k-2}+1}\circ
T^h,...,f_{2^k -1}.f_{2^k -1}\circ T^h)(x)\bigg|\bigg)\\
&\leq C.\bigg(\frac{1}{H} + \frac{1}{H}\sum_{h=1}^H
\frac{1}{N^{k-2}}\sum_{i_1,..., i_{k-2}=
0}^{N-1}\\
&\bigg|\frac{1}{N}\sum_{i_k=0}^{N-1}A_{N,(i_1,i_2,...,i_{k-2},i_k)}(f_{3.2^{k-2}+1}.f_{3.2^{k-2}+1}\circ
T^h,...,f_{2^k -1}.f_{2^k -1}\circ T^h)(x)\bigg|\bigg)\\
&\leq C.\bigg(\frac{1}{H} + \bigg(\frac{1}{H}\sum_{h=1}^H
\frac{1}{N^{k-2}}\sum_{i_1,..., i_{k-2}=
0}^{N-1}\\
&\bigg|\frac{1}{N}\sum_{i_k=0}^{N-1}A_{N,(i_1,i_2,...,i_{k-2},i_k)}(f_{3.2^{k-2}+1}.f_{3.2^{k-2}+1}\circ
T^h,...,f_{2^k -1}.f_{2^k -1}\circ
T^h)(x)\bigg|^2\bigg)^{1/2}\bigg)
\end{aligned}
\]
So by the induction assumption we have
\[
\begin{aligned}
&\limsup_{N}\frac{1}{N^{k-2}}\sum_{i_1,..., i_{k-2}=
0}^{N-1}\sup_{t}\bigg|\frac{1}{N}\sum_{i_k=0}^{N-1}A_{N,(i_1,i_2,...,i_{k-2},i_k)}(f_{3.2^{k-2}+1},...,f_{2^k
-1})(x)e^{2\pi ii_kt}\bigg|^2 \\
&\leq C.\bigg(\frac{1}{H} + \bigg(\frac{1}{H}\sum_{h=1}^H
\limsup_N\frac{1}{N^{k-2}}\sum_{i_1,..., i_{k-2}=
0}^{N-1}\\
&\bigg|\frac{1}{N}\sum_{i_k=1}^{N-1}A_{N,(i_1,i_2,...,i_{k-2},i_k)}(f_{3.2^{k-2}+1}.f_{3.2^{k-2}+1}\circ
T^h,...,f_{2^k -1}.f_{2^k -1}\circ
T^h)(x)\bigg|^2\bigg)^{1/2}\bigg)\\
&\leq C .\bigg(\frac{1}{H} + \bigg(\frac{1}{H}\sum_{h=1}^H
Min_{\{3.2^{k-2}+1 \leq j \leq 2^{k}-1\}}\||f_jf_j\circ
T^h|\|_{k-2}^2\bigg)^{1/2}\bigg)
\end{aligned}
\]
By using the monotonicity in $\alpha$ of the fractions
$\big(\frac{1}{P}\sum_{h=1}^H |u_h|^{\alpha}\big)^{1/\alpha},$ we
have
\[
\begin{aligned}
&\limsup_{N}\frac{1}{N^{k-2}}\sum_{i_1,..., i_{k-2}=
0}^{N-1}\sup_{t}\bigg|\frac{1}{N}\sum_{i_k=0}^{N-1}A_{N,(i_1,i_2,...,i_{k-2},i_k)}(f_{3.2^{k-2}+1},...,f_{2^k
-1})(x)e^{2\pi ii_kt}\bigg|^2 \\
&\leq C .\bigg(\frac{1}{H} + \bigg(\frac{1}{H}\sum_{h=1}^H
Min_{\{3.2^{k-2}+1 \leq j \leq 2^{k}-1\}}\||f_jf_j\circ
T^h|\|_{k-2}^{2^{k-2}}\bigg)^{1/2^{k-2}}\bigg)
\end{aligned}
\]
By taking now the $\limsup_H$ of the last term we get
\begin{equation}
\begin{aligned}
&\limsup_{N}\frac{1}{N^{k-2}}\sum_{i_1,..., i_{k-2}=
0}^{N-1}\sup_{t}\bigg|\frac{1}{N}\sum_{i_k=0}^{N-1}A_{N,(i_1,i_2,...,i_{k-2},i_k)}(f_{3.2^{k-2}+1},...,f_{2^k
-1})(x)e^{2\pi ii_kt}\bigg|^2 \\
&\leq C.Min_{\{3.2^{k-2}+1 \leq j \leq 2^{k}-1\}}\||f_j|\|_{k-1}^2
\end{aligned}
\end{equation}
Thus if one of the functions $f_j$ belongs to
$\mathcal{Z}_{k-1}^{\perp}$ then the limit in the equation (7) is
equal to zero.
\end{proof}

\smallskip

\noindent{\bf End of the proof of theorem 1}
\smallskip

 We just need to finish the induction process the same way we did in \cite{assani} by proving the
 induction assumption for k .
 We consider the averages of $2^k -1$ functions $f_j$, $M_N(f_1, f_2,..., f_{2^k
 -1})(x)$. With the inequality (5) we have
\[\begin{aligned}
&|M_N(f_1, f_2,..., f_{2^k -1})(x)|^2 \\
&\leq \prod_{j=1}^{2^{k-1}-1}\|f_j\|_{\infty}^2
\frac{1}{N^{k-1}}\sum_{i_1,...,
 i_{k-1}=0}^{N-1}\big|\frac{1}{N}\sum_{i_k=0}^{N-1}S_{N,(i_1,i_2,...,i_k)}(f_{2^{k-1}},..., f_{2^k
-1})(x)\big|^2.
\end{aligned}
\]
 By using the same method used to derive the equations (1) and
(3) we get

 \[
 \begin{aligned}
 &\prod_{j=1}^{2^{k-1}-1}\|f_j\|_{\infty}^2
\frac{1}{N^{k-1}}\sum_{i_1,...,
 i_{k-1}=0}^{N-1}\big|\frac{1}{N}\sum_{i_k=0}^{N-1}S_{N,(i_1,i_2,...,i_k)}(f_{2^{k-1}},..., f_{2^k
-1})(x)\big|^2 \\
 &\leq C\frac{1}{N^{k-2}}\sum_{i_2,...,
 i_{k-1}=0}^{N-1}\sup_t\big|\frac{1}{N}\sum_{i_k^{'}=0}^{N-1}A_{N, (i_2,...,i_{k-1},
i_k^{'})}(f_{3.2^{k-2}+1},...,f_{2^k -1})(x)e^{2\pi
i_k^{'}t}\big|^2
\end{aligned}
\]
By using lemma 1 and (8) one concludes that
\[\begin{aligned}
&\limsup_N\frac{1}{N^{k-1}}\sum_{i_1,...,
 i_{k-1}=0}^{N-1}\big|\frac{1}{N}\sum_{i_k=0}^{N-1}S_{N,(i_1,i_2,...,i_k)}(f_{2^{k-1}},..., f_{2^k
-1})(x)\big|^2\\
&\leq C\frac{1}{N^{k-2}}\sum_{i_2,...,
 i_{k-1}=0}^{N-1}\sup_t\big|\frac{1}{N}\sum_{i_k^{'}=0}^{N-1}A_{N, (i_2,...,i_{k-1},
i_k^{'})}(f_{3.2^{k-2}+1},...,f_{2^k -1})(x)e^{2\pi
i_k^{'}t}\big|^2\\
 &\leq C.Min_{\{3.2^{k-2}+1 \leq j \leq
2^{k}-1\}}\||f_j|\|_{k-1}^2
\end{aligned}
\]
 By symmetry on the indices $i_1, i_2, ..., i_k$ one obtains the
 following inequality for the $2^{k-1}$ functions $f_j$
 \[
 \begin{aligned}
&\limsup_N\frac{1}{N^{k-1}}\sum_{i_1,...,
 i_{k-1}=0}^{N-1}\big|\frac{1}{N}\sum_{i_k=0}^{N-1}S_{N,(i_1,i_2,...,i_k)}(f_{2^{k-1}},..., f_{2^k
-1})(x)\big|^2\\
&\leq C Min_{\{2^{k-1}\leq j\leq 2^k -1\}}\||f_j|\|_{k-1}^2
\end{aligned}
\]
By applying this last inequality to any set of $2^{k-1}$ functions
functions $g_j$ that we can label from $3.2^{k-1} +1$ to
$2^{k+1}-1$ instead of $1$ to $2^k -1$ we obtain our induction
assumption for k.
\smallskip
Thus the averages $M_N(f_1, f_2,..., f_{2^k -1})(x)$ converge a.e.
to zero if one of the functions $f_j\in Z_{k-1}^{\perp}$ (using
the symmetry of the indices). Combining this result with the
pointwise convergence when all functions are in $Z_{k-1}$
mentioned in \cite{Host-Kra}, (see \cite{Leibman}) this ends the
proof of theorem 1.

\section{proof of theorem 2}
\smallskip

 The proof of theorem 2 follows from the path we used. We showed
 that if one of the functions $f_j$ is in the orthocomplement of
 the $Z_{k-1}$ factor then the averages of these $2^k -1$
 functions converge a.e to zero. Thus the limit is given by the
 pointwise convergence when all functions are in the factor
 $Z_{k-1}$.


\begin{thebibliography}{99}
\bibitem{assani}{\bf{I. Assani}}:``Pointwise convergence of
averages along cubes'', preprint.

\bibitem{Host-Kra}{\bf{B. Host and B. Kra}}: ``Nonconventional
ergodic averages and nilmanifolds'', to appear in Annals of Math.

\bibitem{Bergelson}{\bf{V. Bergelson}}:``The multifarious Poincare
Recurrence theorem,'' \textit{Descriptive Set Theory and Dynamical
Systems,}, Eds M. Foreman, A.S. Kechris, A. Louveau, B. Weiss.
Cambridge University Press, New York (2000), 31-57.

\bibitem{Furstenberg-Weiss}{\bf H. Furstenberg and B. Weiss}:``A
mean ergodic theorem for $\dis \frac{1}{N}\sum_{n=1}^N
f(T^nx)g(T^{n^2}x)$'',\textit{Convergence in Ergodic Theory and
Probability}, Eds: Bergelson/March/Rosenblatt, Walter de
Gruyter-Co, Berin, New York (1996), 193-227.

\bibitem{Khintchine}{\bf{A. Y. Khintchine}}:``Eine Verscharfung des
Poincareschen "Wiederkehrsatzes",''\textit{Comp. Math.}, 1,
(1934), 177-179.

 \bibitem{Kuipers-Niederreiter}{\bf{L. Kuipers and H. Niederreiter}}:
  \textit{Uniform Distribution of Sequences.} John Wiley \& Sons, 1974.

 \bibitem{Leibman}{\bf {A. Leibman}}:
\textit{Pointwise convergence of ergodic averages for polynomial
sequences of rotations of a nilmanifold}:
http://www.math.ohio-state.edu/~leibman/preprints/
  \end{thebibliography}
  \end{document}